\documentclass{amsart}

\usepackage{enumitem}   

\setlist[enumerate]{itemsep=.2em,topsep=.2em,leftmargin=1.25em,itemindent=2.0em}

\usepackage{amssymb}
\usepackage[all]{xy}
\usepackage{hyperref}

\newtheorem{thm}{Theorem}

\newtheorem{lem}[thm]{Lemma}
\newtheorem{cor}[thm]{Corollary}

\newtheorem{prop}[thm]{Proposition}

\theoremstyle{definition}

\newtheorem{say}[thm]{}
\newtheorem{exmp}[thm]{Example}

\newtheorem{defn-thm}[thm]{Definition--Theorem}  
\newtheorem{defn-lem}[thm]{Definition--Lemma}  

\theoremstyle{remark}

\newcommand{\n}[0]{{\mathbb N}}

\newcommand{\cc}[0]{{\mathbb C}}
\newcommand{\dd}[0]{{\mathbb D}}
\newcommand{\bdd}[0]{\overline{\mathbb D}}

\newcommand{\p}[0]{{\mathbb P}}

\def\into{\DOTSB\lhook\joinrel\to}

\def\loccoh#1.#2.#3.#4.{H^{#1}_{#2}(#3,#4)}

\DeclareMathAlphabet{\mathchanc}{OT1}{pzc}%
                                {m}{it}

\usepackage[all]{xy}\xyoption{dvips}

\newcommand{\dist}[0]{\operatorname{dist}}

\begin{document}
\bibliographystyle{amsalpha}

   \title[Comment]{Comment on:  Dense entire curves ...\\  by F.~Campana and J.~Winkelmann}
   \author{J\'anos Koll\'ar}

\maketitle

We explain how  modifications of the methods of \cite{MR2150882, c-w-dense} give the following.

\begin{thm} \label{nev.hilb.thm}
  Let $X$ be a smooth, proper variety over $\cc$. Then there is a holomorphic map $\phi:\dd\to X$ such that, for every finite cover   $p:Y\to X$, the map $\phi$   lifts to  $\phi_Y:\dd\to Y$ iff $p$ is \'etale.

  If $X$ is rationally connected, then there is also such a map
   $\phi:\cc\to X$.
\end{thm}

By a {\it finite cover} we mean a dominant, finite morphism from an irreducible variety $Y$ to $X$; thus ramification is allowed.  If  $X$ is  rationally connected, then it is simply connected, hence a  nontrivial cover must have ramification.

The open (resp.\ closed) complex unit disc is denoted by $\dd $ (resp.\ $\bdd$).
A map $g:\bdd \to X$ is holomorphic if it extends holomorphically to a  larger disc.

The  $\phi:\cc\to X$ variant of Theorem~\ref{nev.hilb.thm} is called the {\it Nevanlinna version of the Hilbert property} in  \cite{c-w-dense}, following \cite[\S 2.4]{MR3648511}.
It is proved in \cite[9.4]{c-w-dense} for Galois covers.  Also,
\cite[9.9]{c-w-dense} shows that for every  nontrivial $p:Y\to X$ there is a
$\phi_Y:\cc\to X$ that can not be lifted to $Y$, but the  $\phi_Y$ constructed there does depend on $Y$. See Example~\ref{galois.dexmp} for the difference between Galois covers and arbitrary covers in this context.

\subsection*{Background results}{\ }

The first ingredients of the proof are
the following. They are not stated, but essentially proved in \cite{MR2150882} and \cite[Sec.5]{c-w-dense}; see  especially \cite[5.1]{c-w-dense}.

\begin{prop} \label{dd.prop.cw}
  Let $X$ be a smooth, proper variety over $\cc$ and fix any Riemannian metric on $X$.
  Let $\psi_i:\bdd\to X$ be countably many holomorphic maps 
  and $\epsilon_i>0$. Then there is  a holomorphic map
   $\phi:\dd\to X$ and  embeddings $\tau_i:\bdd\into \dd$ such that
  $$
  \dist\bigl(\phi(\tau_i(z)), \psi_i(z)\bigr)<\epsilon_i
  $$
  for every $z\in \bdd$ and every $i$. \qed
\end{prop}

In the rationally connected case, we get a similar results for discs that come from rational curves.

\begin{prop} \label{prop.cw}
  Let $X$ be a smooth, proper, rationally connected variety over $\cc$ and fix any Riemannian metric on $X$.
  Let $g_i:\p^1\to X$ be countably many free rational curves, $\sigma_i:\bdd\into \p^1$ closed discs
  and $\epsilon_i>0$. Set $\psi_i:=g_i\circ \sigma_i$. Then there is  a holomorphic map
   $\phi:\cc\to X$ and  embeddings $\tau_i:\bdd\into \cc$ such that
  $$
  \dist\bigl(\phi(\tau_i(z)), \psi_i(z)\bigr)<\epsilon_i
  $$
  for every $z\in \bdd$ and every $i$. \qed
\end{prop}

We also use the following  Lefschetz-type result for surjectivity of 
$\pi_1$.

\begin{prop}\cite{MR820315, MR2011744}\label{prop.k-rcfg}
  Let $X$ be a smooth, proper variety over $\cc$.
  Then there is a family of smooth  curves
$ W\leftarrow C_W\stackrel{u}{\to} X$
  with the following property.

 For every finite  cover 
 $p:Y\to  X$ there is a dense, open $W_Y\subset W$
 such that for every $w\in W_Y$
 \begin{enumerate}
 \item the fiber product $C_w\times_X Y$ is irreducible, smooth, and
 \item $C_w\times_X Y\to C_w$ is \'etale  iff $p$ is \'etale.
   \item Furthermore,  if $X$ is rationally connected, we can choose
  $C_W=W\times \p^1$.
   \end{enumerate}
 \end{prop}

Proof. If $X$ is projective, then
the universal family of smooth curve sections with linear spaces works by
\cite{MR820315}. In the proper case we use the images of linear space sections of a projective modification $X'\to X$.

The rationally connected case is proved in \cite[Cor.7]{MR2011744}.\qed
\medskip

We also need some lemmas  that countably many maps suffice  in Theorem~\ref{nev.hilb.thm} if $\dim X=1$.

\begin{lem}\label{RS.uc.lem}
  Let $S$ be a  Riemann surface with universal cover
  $\pi:\dd\to S$. Let $p:S'\to S$ be a finite cover such that $\pi$ lifts to
  $\pi':\dd\to S'$. Then $p$ is unramified.
\end{lem}

Proof. After base change to $\dd$ we have
$\tilde\pi: \dd\times_S\dd\to \dd\times_SS'$.
Here $\dd\times_S\dd $  is a disjoint union of countably many copies of $\dd$, and each one gives a section of the coordinate projection $\tilde p: \dd\times_SS'\to \dd$.
Since the group of deck transformations acts transitively on the connected components of $\dd\times_SS'$, we see that $\tilde p$ is unramified, and so is $p$. \qed

\begin{cor} \label{RS.uc.cor}  Using the notation of Lemma~\ref{RS.uc.lem}, 
     let $\pi_m:\bdd(1-1/m)\to S$ denote the restriction of 
    $\pi$ for $m\in\n$.   If infinitely many of the $\pi_m$ lift to
    $\pi'_m:\bdd(1-1/m)\to S'$, then $p$ is unramified.
\end{cor}

Proof. Fix $z_0\in \dd$ such that $p$ is unramified over   $\pi(z_0)$.
Then $\pi'_m(z_0)$ is the same point $z'_0\in S'$ for infinitely many $m$.
These define a lifting $\pi':\dd\to S'$. So $p$ is unramified by Lemma~\ref{RS.uc.lem}. \qed

\begin{lem} \label{P1.uc.cor}  There are countably many discs
  $\sigma_j:\bdd\into \cc\p^1$ such that a  cover  $p:S'\to \cc\p^1$ is trivial iff every $\sigma_j$ lifts to $\sigma'_j:\bdd\into  S'$.
\end{lem}

Proof. Choose  discs  $\sigma_j$ such that every finite subset of $\cc\p^1$ is contained in the image of some $\sigma_j(\dd)$. For example, we can use the interiors and exteriors of all circles $|z|=r$ where either $r$ or $r^{-1}$ is a natural number.

For any  $p:S'\to \cc\p^1$ there is a   $\sigma_j$ whose image contains all branch points. Then  $p^{-1}(\sigma_j(\dd))$ is irreducible. So  $\sigma_j$ lifts iff
$p$ is an isomorphism over $\sigma_j(\dd)$, hence everywhere. \qed

\begin{say}[Proof of Theorem~\ref{nev.hilb.thm}] \label{nev.hilb.thm.pf} We choose curves 
  $C_{w_*}\to X$, discs  $\tau_*:\bdd\to C_{w_*}$  and $\epsilon_*>0$ as follows.
\begin{enumerate}
\item First we pick a  dense set  of points $w_i\in W$.  (Zariski dense is sufficient.)
\item  Next we choose  discs  $\tau_{ij}:\bdd\to C_{w_i}$ as in
  Corollary~\ref{RS.uc.cor}  (resp.\ Lemma~\ref{P1.uc.cor}).  
\item   Set $\epsilon_k=2^{-k}$.
\end{enumerate}

We apply Propositions~\ref{dd.prop.cw}--\ref{prop.cw} to the triply infinite family of all
$$
\tau_{ijk}:=\tau_{ij},\quad \epsilon_{ijk}:=\epsilon_k.
$$
We claim that any resulting  $\phi:\dd\to X$  (resp.\ $\phi:\cc\to X$) works for Theorem~\ref{nev.hilb.thm}.

To see this, choose any finite, non-\'etale  cover $p:Y\to X$.
By Proposition~\ref{prop.k-rcfg}  and (\ref{nev.hilb.thm.pf}.1)
there is a $w_i$ such that   
$C_{w_i}\times_XY$ is smooth, irreducible  and
$$
C_{w_i}\times_XY\to C_{w_i}
$$
is non-\'etale.  

Next, by  (\ref{nev.hilb.thm.pf}.2) there is a
$\tau_{ij}:\bdd\to C_{w_i}$ that does not lift to $C_{w_i}\times_XY$.

  If $\phi$ lifts to $\phi_Y:\dd\to Y$ (resp.\ $\phi_Y:\cc\to Y$) then, fixing $i, j$ and letting $k\to \infty$,   a subsequence of the holomorphic maps
  $$
  \phi_Y\circ \tau_{ijk}:\bdd\to Y
  $$
  converges to a lifting of $\tau_{ij}:\bdd\to X$.
  Thus we get a lifting of $\tau_{ij}$ to  $C_{w_i}\times_XY$.
  This is impossible by Corollary~\ref{RS.uc.cor}  (resp.\ Lemma~\ref{P1.uc.cor}).  \qed
\end{say}

The difference between lifting to Galois covers and arbitrary covers is illustrated by the following.

\begin{exmp}\label{galois.dexmp}
  Let $\phi:\dd\to S$ be a holomorphic map to a Riemann surface that is not simply connected.
\begin{enumerate}
\item  If $\phi$ is surjective and \'etale, then it can not be  lifted to any non-\'etale, Galois cover $S'\to S$.
\item
  If $\phi^{-1}(s)$ is finite for some $s\in S$, then $\phi$ can be lifted to some  non-\'etale, non-Galois cover $S'\to S$.
\end{enumerate}
For example let $u:\dd\to S$ be the universal cover and consider the smaller discs
$\bdd(1-c)\subset \dd$. We get
$\phi_c: \bdd(1-c)\to S$.  Then  every $\phi_c$ can be lifted to some  non-\'etale, non-Galois cover $S'\to S$, but if $c\ll 1$  then   $\phi_c$ can not be  lifted to any non-\'etale, Galois cover $S'\to S$.

Claim (2) shows that a feature of the proof of Theorem~\ref{nev.hilb.thm}---that the image of $\phi$ returns to some neighborhoods infinitely many times---is necessary. 
  \end{exmp}


\def\cprime{$'$} \def\cprime{$'$} \def\cprime{$'$} \def\cprime{$'$}
  \def\cprime{$'$} \def\dbar{\leavevmode\hbox to 0pt{\hskip.2ex
  \accent"16\hss}d} \def\cprime{$'$} \def\cprime{$'$}
  \def\polhk#1{\setbox0=\hbox{#1}{\ooalign{\hidewidth
  \lower1.5ex\hbox{`}\hidewidth\crcr\unhbox0}}} \def\cprime{$'$}
  \def\cprime{$'$} \def\cprime{$'$} \def\cprime{$'$}
  \def\polhk#1{\setbox0=\hbox{#1}{\ooalign{\hidewidth
  \lower1.5ex\hbox{`}\hidewidth\crcr\unhbox0}}} \def\cdprime{$''$}
  \def\cprime{$'$} \def\cprime{$'$} \def\cprime{$'$} \def\cprime{$'$}
\providecommand{\bysame}{\leavevmode\hbox to3em{\hrulefill}\thinspace}
\providecommand{\MR}{\relax\ifhmode\unskip\space\fi MR }
\providecommand{\MRhref}[2]{%
  \href{http://www.ams.org/mathscinet-getitem?mr=#1}{#2}
}
\providecommand{\href}[2]{#2}

\bigskip

  Princeton University, Princeton NJ 08544-1000, 

\email{kollar@math.princeton.edu}

\end{document}